\documentclass[12pt]{amsart}
\usepackage{amsfonts}
\usepackage{bbm}
\usepackage{mathrsfs}
\usepackage{amssymb,amsmath,amsfonts,amsthm}

\numberwithin{equation}{section}

\newcommand{\qe}{\hfill\Box}
\newtheorem{thm}{\hskip\parindent {Theorem}}[section]

\newtheorem{defn}{\hskip\parindent {Definition}}[section]

\def\rr{\Bbb{R}}
\def\cc{\Bbb{C}}

\def\zz{\Bbb{Z}}

\def\M{\mathcal{M}}

\def\S{\mathcal{S}}
\def\R{\mathcal{R}}

\def\rrnm{\rr^n\times \rr^m}
\numberwithin{equation}{section}

\begin{document}

\baselineskip 17.2pt \hfuzz=6pt

\title[ An atomic decomposition characterization of $H^p_F({\rr}^{n}\times{\rr}^m)$]
{An atomic decomposition characterization of flag Hardy spaces
$H^p_F({\rr}^{n}\times{\rr}^m)$ with applications}

\bigskip

\author{ Xinfeng Wu   }

\address{Xinfeng Wu \\
    Department of Mathematics\\
    China University of Mining and Technology (Beijing) \\Beijing 100083,
 China}
     \email{wuxf@cumtb.edu.cn}

\thanks{
This research is supported by NNSF-China No. 11101423, 11171345.}


\begin{abstract}
In this paper, we give an atomic decomposition characterization of
flag Hardy spaces $H^p_F({\rr}^n\times {\rr}^m)$ for $0<p\le 1$,
which were introduced in \cite{hl1}. A remarkable feature of atoms of such flag Hardy spaces is that these atoms have only
partial cancellation conditions. As an application, we prove a boundedness criterion for operators on flag Hardy spaces.

\end{abstract}

\maketitle

\section{Introduction and the main result}
Atomic decomposition is a significant tool in studying various
function spaces and operators arising in harmonic analysis and
wavelet analysis (see see Meyer \cite{Y} and Coifman-Meyer
\cite{cm}, etc.). The atoms can be viewed as building blocks of
functions spaces and the problem of boundedness on function spaces
can simply be reduced to the uniform boundedness of operators on
such building blocks ``atoms''.

The atomic decompositions for Hardy spaces in one dimension were
first constructed by Coifman in 1974 (see \cite{Co}) and later were
extended to higher dimensions by Latter \cite{La}. The
multi-parameter case is more complicated. Indeed, a product atom  is
no longer supported in a rectangle (the direct analogue of the
cube), but in a general open set. The celebrated $(p,2)$-atomic
decomposition of product Hardy spaces $H^p({\rr}\times {\rr})$ were
constructed by Chang-R. Fefferman in \cite{cf,cf2,cf3}. Recently,
the general $(p,q)$ atomic decomposition were given by Han, Lu and
Zhao in \cite{hlz}.

The flag singular integral operators were first introduced by M\"{u}ller, Ricci and Stein
when they studied the Marcinkiewicz multiplier on the Heisenberg groups in \cite{mrs}.
To study the $\square_b$-complex on certain   CR submanifolds
of ${\cc}^n$, in 2001, Nagel, Ricci and Stein \cite{nrs} studied a
class of product singular integrals with flag kernel. They proved, among other things, the $L^p$ boundedness of flag singular integrals.
More recently, Nagel-Ricci-Stein-Wainger in \cite{nrsw1,nrsw2} have generalized these results to a more general setting, namely, homogeneous group. For other related results, see \cite{g1,g2}.

For $0<p\leq 1$, Han and Lu \cite{hl1} developed a unified approach of Hardy spaces
with respect to flag multi-parameter structure. The $H^p$ and
$H^p-L^p$ boundedness of flag singular integral operators and the duality of $H^p$ were established via
the discrete Littlewood-Paley-Stein analysis and discrete Calder\'{o}n's identity. These ideas were later carried out to the weighted case in \cite{dhlw}. However, an atomic decomposition characterization for the multiparameter flag Hardy spaces is still absent.

The purpose of the present paper is to establish an atomic decomposition
characterization for flag Hardy spaces $H^p_F({\rr}^n\times{\rr}^m)$
introduced in \cite{hl1}. As pointed out in
\cite{hl1}, the main difficulty lies in the fact that the underling
multiparameter structure in the flag case is {\it not explicit}, but
{\it implicit}. As a consequence, the classical cancellation conditions on atoms are not available for the flag
structure. To overcome this difficulty, we lift the flag structure to product
structure in higher dimensions. More precisely, two types of lifting must be considered separately according to different geometric shapes of the rectangles. These two different kinds of lifting lead naturally to two different cancellation conditions.

Before stating our main result, we first recall some basic notions
and notations in \cite{hl1}.

The product test function class ${\S}_{\infty}({\rr}^{n+m}\times
{\rr}^m)$ is the collection of all functions $f\in
{\S}({\rr}^{n+m}\times {\rr}^m)$ with
$$
\int_{{\rr}^{n+m}} f(x,y,z)x^{\alpha}y^{\beta}dxdy=\int_{{\rr}^m}
f(x,y,z)z^{\gamma}dz=0
$$
for all multi-indices $\alpha,\beta,\gamma$ of nonnegative integers.

\begin{defn}\label{def flag test function:f}
A function $f(x,y)$ defined on ${\rr}^n\times {\rr}^m$ is said to be
a test function in ${\S}_F({\rr}^n\times{\rr}^m)$ if there exists a
function $f^{\#}\in {\S}_{\infty}({\rr}^{m+n}\times {\rr}^m)$ such
that
\begin{eqnarray}\label{eq: test function:f}
f(x,y)=\int_{{\rr}^m}f^{\#}(x,y-z,z)dz,
\end{eqnarray}
the norm is defined by
$$
\|f\|_{{\S}_F({\rr}^{n}\times
{\rr}^m)}=\inf\{\|f^{\#}\|_{{\S}_{\infty}({\rr}^{n+m}\times{\rr}^m)}:
\mbox{for all representations of} \ f\ {\rm in}\ (\ref{eq: test
function:f})\}.
$$
Denote by $({\S}_F)'$ the dual of ${\S}_F$.
\end{defn}

Let $\psi^{(1)}\in\S({\rr}^{n+m})$ ,$\psi^{( 2)}\in \S({\rr}^{m})$
be Schwartz functions supported in the unit balls of ${\rr}^{n+m}$
and ${\rr}^m$  and satisfy
$$
\sum_j|{\widehat{\psi^{(1)}}(2^{-j}\xi_1,2^{-j}\xi_2)}|^2=1
$$
for all $(\xi_1,\xi_2)\in {\rr}^n\times {\rr}^m\backslash
\{(0,0)\}$, and
$$
\sum_k|\widehat{\psi^{(2)}}(2^{-k}\eta)|^2=1
$$
for all $\eta\in {\rr}^m\backslash\{0\}$ and the moment conditions
$$
\int \psi^{(1)}(x,y)x^{\alpha}y^{\beta}dxdy=0=\int
\psi^{(2)}(z)z^{\gamma}dz
$$
for multi-indices $|\alpha|,|\beta|,|\gamma|\le M_0$, where $M_0$ is
a large integer satisfying $M_0\ge [2/p-1/2]\max\{m,n\}$. Set $\psi_{j,k}=\psi_j^{(1)}*_2\psi^{(2)}_k,$
where $*_2$ denotes the convolution in the second variable of
$\psi_j^{(1)}$ with $\psi^{(2)}_k$, where $\psi_j(x_1,x_2)=2^{-j(n+m)}\psi(2^{-j}x_1,2^{-j}x_2)$ and $\psi_k(x_3)=2^{-km}\psi(2^{-k}x_3)$.

For $f\in ({\S}_F)'({\rr}^n\times {\rr}^m)$, the Littlewood-Paley-Stein
square function of flag type $g_{F}(f)$ is defined by
$$
g_{F}(f)(x,y)=\left(\sum_{j\in\zz}\sum_{k\in\zz}|\psi_{j,k}*f(x,y)|^2\right)^{1/2}.
$$
The multiparameter flag Hardy spaces can be defined as follows.
\begin{defn}\label{def flag hardy:f}
For $0<p\le 1$, define flag Hardy spaces $H_F^p({\rr}^n\times {\rr}^m)$ by
$$H_F^p({\rr}^n\times {\rr}^m)=\{f\in ({\S}_F)':g_{F}(f)\in
L^p({\rr}^n\times {\rr}^m)\}$$ and the $H^p_F$ norm by
$$
\|f\|_{H^p_F({\rr}^n\times {\rr}^m)}=\|g_{F}(f)\|_{L^p({\rr}^n\times {\rr}^m)}.
$$
\end{defn}

Now we give the definition of multiparameter flag $H^p_F({\rr}^{n}\times {\rr}^m)$
atoms.

\begin{defn}\label{def flag atom:f}
Suppose $0<p\le1$. A multiparameter flag $H^p$-atom (or $H^p_F$-atom) is a function
$a(x_1,x_2)$ on ${\rr}^n\times {\rr}^m$ supported in some open set
$\Omega$ of finite measure such that
\begin{enumerate}
\item $\|a\|_{L^2({\rr}^n\times{\rr}^m)}\le |\Omega|^{1/p-1/2}$;
\item $a$ can further be decomposed into ``flag particles'' $a_R$
as follows.
$$a=\sum_{R\in \Omega}a_R,$$
where $a_R$ is supported in $6R=6(I_R\times J_R)$ for some $R\subset
\Omega$ with $l(I_R)=2^{j_R}$ and $l(J_R)=2^{k_R}$; the flag
particles can be lifted to $a^{\#}_R$ and $\widetilde{a^{\#}_R}$ by
\begin{align*}
a_R(x_1,x_2)=&\int a_R^{\#}(x_1,x_3;x_2-x_3)dx_3=\int \widetilde{a_R^{\#}}(x_1,x_2-x_3;x_3)dx_3
\end{align*}
and $a^{\#}_R(x_1,x_3;x_2)$, $\widetilde{a^{\#}_R}(x_1,x_2;x_3)$ are
``pseudo-product particles'' which satisfy the following properties.
If $2^{k_R}\ge 2^{j_R}$, then
\begin{enumerate}
\item[\rm(a)] {\rm supp} $a^{\#}_R\subset R^{\#}:=I_R\times
\hat{I}_R\times J_R$, where $\hat{I}_R$ is a cube in ${\rr}^m$
centered at origin with the same sidelength as $I_R$;
\item[\rm(b)] $a_R^{\#}$ satisfies the following moment
conditions:
\begin{align*}
\iint_{{\rr}^m\times{\rr}^n} a_R^{\#}(x_1,x_3;x_2)x_1^{\alpha}
x_3^{\beta}dx_1dx_3&=0, {\rm for\ all} \ x_2 \in {\rr}^m {\rm and}\ |\alpha|,\, |\beta|\le
k_pn-1,
\end{align*}
and
$$\int_{{\rr}^m}a_R^{\#}(x_1,x_3;x_2)x_2^{\gamma}dx_2=0,\quad
{\rm for\ all}\ (x_1,x_3)\in {\rr}^{n+m}\ {\rm and}\ |\gamma|\le
k_pm-1,$$ where $k_p\le [2/p-1/2]$. \item[\rm(c)]$a_R^{\#}$ is
$C^{k_pn}$ in $x_1$ and $x_3$ and $C^{k_pm}$ in $x_2$ with
$\|a_R^{\#}\|_{\infty}\le d_R$,
$$
\left\|\frac{\partial^{kn} a_R^{\#}}{\partial
x_1^{kn}}\right\|_{\infty},\quad \left\|\frac{\partial^{kn}
a_R^{\#}}{\partial x_3^{kn}}\right\|_{\infty}\ \le
\frac{d_R}{|I_R|^k},\
\mbox{and}\
\left\|\frac{\partial^{km} a_R^{\#}}{\partial x_2^{km}}\right\|_{\infty}\
\le \frac{d_R}{|J_R|^k}$$ for each $k\le k_p$ and
\begin{eqnarray}\label{pseudo particle:f}
\sum d_R^2\, |R|\,|\hat{I}_R|^2 \le A|\Omega|^{1-2/p};
\end{eqnarray}
\end{enumerate}
If $2^{j_R}\ge 2^{k_R}$, then {\rm supp}
$\widetilde{a^{\#}_R}\subset \widetilde{R^{\#}}:=I_R\times
\tilde{I}_R\times J_R$, where $\tilde{I}_R:=\hat{I_R}+c_{J_R}$
($c_{J_R}$ denotes the center of $j_R$) and it satisfies similar
properties as {\rm (b)} and {\rm (c)} except for the role of $x_2$
and $x_3$ exchanged.
\end{enumerate}
\end{defn}

Our main result is
the following

\begin{thm}\label{main flag atom:f}
Let $0<p\le 1$. If $f\in H^p_F({\rr}^n\times{\rr}^m)$ if and only if
$f$ can be written as $f=\sum \lambda_k a_k$ where $a_k$ are flag
$p$-atoms and the series converges in the norm of $H^p_F(\rrnm).$ Moreover,
$$\|f\|_{H^p_F}\approx \inf\lbrace (\sum\lambda_k^p)^{\frac{1}{p}}:\ \mbox{for all}\ f=\sum \lambda_k a_k\rbrace.$$
In particular, if $f$ is in both $L^2(\rrnm)$ and $H^p_F(\rrnm),$ then $f=\sum \lambda_k a_k $ where the series converges in both $L^2(\rrnm)$ and $H^p_F(\rrnm),$ and $(\sum\lambda_k^p)^{\frac{1}{p}}\lesssim \|f\|_{H^p_F}.$
\end{thm}

As an application of Theorem \ref{main flag atom:f}, we obtain the following boundedness criterion for operators on flag Hardy spaces.
\begin{thm}\label{main:application}
Let $T$ be a $L^2$ bounded linear operator. Then $T$ is bounded from
$H^p_F(\rrnm)\ 0<p\le 1$ to $L^p(\rrnm)\ 0<p\le 1$ if and only if
$\|Ta\|_{L^p}\le C$ uniformly for all $H^p_F(\rrnm)$ atoms $a,$ and $T$ is bounded on
$H^p_F(\rrnm)\ 0<p\le 1$ if and only if
$\|Ta\|_{H^p_F}\le C$ uniformly for all $H^p_F(\rrnm)$ atoms $a.$
\end{thm}

\section{Proof of Theorem \ref{main flag atom:f}}
First we show the
necessity by constructing an atomic decomposition. Suppose $f\in H^p_F(\rrnm).$ We need the
maximal square function defined by
\begin{align*}
g_F^{sup}(f)(x,y)&:=\left\{\sum_{j,k}\sum_{I,J}\sup_{u\in I,v\in
J}|\psi_{j,k}*f(u,v)|^2\chi_I(x)\chi_J(y)\right\}^{1/2},
\end{align*}
where $I\subset {\rr}^n,J\subset{\rr}^m$ are dyadic cubes with
side-length $l(I)=2^{-j}$ and $l(J)=2^{-k}+2^{-j}$, $\chi_I,\chi_J$
are the indicator functions of $I$ and $J$, respectively. It has been shown in \cite{hl1} that
\begin{align}\label{eq: small square function:f}
\|g_F^{sup}(f)\|_p\simeq\|f\|_{H^p_F}.
\end{align}

Next, set
$$
\Omega_i=\{(x,y)\in \rr^n\times \rr^m:g_F^{sup}(f)(x,y)> 2^i\}
$$
and
$$
\widetilde{\Omega}_i=\{(x,y)\in \rr^n\times
\rr^m:\M_s(\chi_{\Omega_i})(x,y)> \frac{1}{100}\},
$$
where $\M_s$ is the strong maximal operator. Obviously, $
\Omega_i\subset\widetilde{\Omega}_i$. By the $L^2$ boundedness of
$\M_s$,
\begin{align}\label{measure of Omega:f}
|\widetilde{\Omega}_i|\leq C |\Omega_i|.
\end{align}
Denote
$$
\mathcal{R}_i=\left\{R=I\times J:|R\cap \Omega_i|>\frac{1}{2}|R|,\
|R\cap \Omega_{i+1}|\leq \frac{1}{2}|R|\right\}
$$
where $I,J$ are dyadic cubes in ${\rr}^n,{\rr}^m$ with side length
$2^{j}$ and $\ 2^{k}+2^j$, respectively. For each $R=I\times J\in
\mathcal{R}_i$, define
\begin{align*}
f_{R}(x,y) = \int_R f_{j,k}(t_1,t_2)\psi_{j,k}(x-t_1,y-t_2)dt_1dt_2,
\end{align*}
where
$
f_{j,k}(t_1,t_2):=f*\psi_{j,k}(t_1,t_2).
$
Set
$$
a_i(x,y)=\sum_{R\in
{\R}_i}a_{R}(x,y)=\frac{C}{2^i|\Omega_i|^{1/p}}\sum_{R\in
{\R}_i}f_{R}(x,y)
$$
and
$
\lambda_i=2^i|\Omega_i|^{1/p}/C,
$
where $C$ is a constant which will be determined later.

From the following Calder\'{o}n representation formula (see
\cite{hl1}),
$
f=\sum_{j,k}f*\psi_{j,k}*\psi_{j,k},
$
where the series converges in $H^p_F(\rrnm)$,
we see that
$$
f=\sum_Rf_R=\sum_i \sum_{R\subset {\R}_i}f_R=\sum_i \lambda_i a_i
$$
It is easy to see that
$$
\sum_i \lambda_i^p=c\sum_i 2^{ip}|\Omega_i|\le C\|g^{\sup}_{F}
(f)\|_{L^p}^p\le C\|f\|_{H^p_F}^p.
$$

Next, we shall show that $a_i$ is a $H^p_F$-atom. Note that for each
$R\in {\R}_i$ and $(y_1,y_2)\in R$,
$
{\rm supp}\, f_{j_R,k_R}\subset {\rm supp}\,\psi_{j_R,k_R}
(\cdot-y_1,\cdot-y_2)\subset 2R\subset \widetilde{\Omega_i}.
$ and
thus $a_i$ is supported in the open set $\widetilde{\Omega}_i$.

By the definition of $g_{F}^{sup}$ and (\ref{measure of Omega:f}),
\begin{eqnarray}\label{eq3: claim:f}
\begin{split}
C2^{2i}|\Omega_i| \ge&
\int_{\tilde{\Omega_i}\backslash\Omega_{i+1}}(g^{sup}_{F})^2(f)(x,y)dx dy\\
\ge&\sum_{R=I\times J\in {\R}_i}(\sup_{u\in I,v\in
J}|\psi_{j,k}*f(u,v)|)^2|(I\times
J)\cap \Omega_i|\\
\ge& \ \frac{1}{2}\sum_{R\in
{\R}_i}|I||J|(\sup_{u\in I, v\in J}|\psi_{j,k}*f(u,v)|)^2,
\end{split}
\end{eqnarray}
where $l(I)=2^j,\ l(J)=2^k+2^j.$

By duality and (\ref{eq3: claim:f}),
\begin{eqnarray}\label{dual}
\begin{split}
\|a_i\|_{2}=&\sup_{\|g\|_2\le 1}\left|\int_{{\rr}^m\times
{\rr}^n}a_i(x,y)g(x,y)dx dy\right|\\
\le& \frac{C}{2^i|\Omega_i|^{1/p}}\sup_{\|g\|_2\le 1}\left|\sum_{R\in{\R}_i}\int_R
f_{j_R,k_R}(t_1,t_2)\widetilde{g_{j_R,k_R}}(t_1,t_2)dt_1dt_2\right|\\
\le& \frac{C}{2^i|\Omega_i|^{1/p}}\left(\sum_{R\in
{\R}_i}|R|(\sup_{u\in I, v\in J}|\psi_{j,k}*f(u,v)|)^2 \right)^{1/2}\\
\le& C|\Omega_i|^{\frac{1}{2}-\frac{1}{p}} \le
C|\widetilde{\Omega}_i|^{\frac{1}{2}-\frac{1}{p}},
\end{split}
\end{eqnarray}
where
$
\widetilde{g_{j_R,k_R}}(t_1,t_2)=
\int_{{\rr}^n\times{\rr}^m}g(x,y)\phi_{j_R,k_R}(x-t_1,y-t_2)dt_1dt_2.
$

Finally, we show that $a_R$ is a ``flag particle''. To this end, write
\begin{align*}
&a_R(x_1,x_2)=\int a_R^{\#}(x_1,x_3;x_2-x_3)dx_3=\int \widetilde{a_R^{\#}}(x_1,x_2-x_3;x_3)dx_3,
\end{align*}
where
\begin{align*}
a_R^{\#}(x_1,x_3;x_2)=\frac{C}{2^i|\Omega_i|^{1/p}}\iint_{R}f_{j,k}(y_1,y_2)
\psi_{j}^{(1)}(x_1-y_1,x_3)\psi^{(2)}_k(x_2-y_2)dy_1dy_2
\end{align*}
and
\begin{align*}
\widetilde{a_R^{\#}}(x_1,x_2;x_3)=\frac{C}{2^i|\Omega_i|^{1/p}}\iint_{R}f_{j,k}(y_1,y_2)
\psi_{j}^{(1)}(x_1-y_1,x_2-y_2)\psi^{(2)}_k(x_3)dy_1dy_2.
\end{align*}
Clearly, $a_R^{\#}$ satisfies the moment conditions in $(x_1,x_3)$ variable
and $x_2$ variable respectively while $\widetilde{a_R^{\#}}$ satisfies
the moment conditions in $(x_1,x_2)$ and $x_3$, respectively.
One can also check that if $k>j$, $a_R^{\#}$ is supported in
$R^{\#}=6(I_R\times \hat{I_R}\times J_R)$ and if $k\le j$,
$\widetilde{a_R^{\#}}$ is supported in
$\widetilde{R^{\#}}:=I_R\times \widetilde{I_R}\times J_R$ where
$$
\hat{I_R}:=\{x_3\in {\rr}^m:|x_3|\le 2^{j}\}\quad {\rm and}\quad
\widetilde{I_R}:=\hat{I_R}+c_{J_R}, \mbox{and}\ c_{J_R}\ \mbox{denotes the center}\ of J_R.
$$

To verify the smoothness conditions for $a_R^{\#}$, we define
$$
d_R:=\frac{A|R|^{1/2}}{|R^{\#}|2^i|\Omega_i|^{1/p}}
\left(\iint_R|f_{j,k}(y_1,y_2)|^2dy_1dy_2\right)^{1/2}.
$$
Then the desired estimates for $a^{\#}_R$ follow from
\begin{gather*}
\|\psi^{(1)}_{j}(x_1-y_1,x_3)\psi^{(2)}_{k}(x_2-y_2)\|_{\infty}\le
\frac{A}{|R^{\#}|},\\
\left\|\left(\frac{\partial^{kn}}{\partial
x_1^{kn}}+\frac{\partial^{kn}}{\partial
x_3^{kn}}\right)\left(\psi^{(1)}_{j}(x_1-y_1,x_3)
\psi^{(2)}_{k}(x_2-y_2)\right)\right\|_{\infty}\le
\frac{A}{|R^{\#}||I_R|^k},\\
\left\|\frac{\partial^{km}}{\partial
x_2^{km}}\left(\psi^{(1)}_{j}(x_1-y_1,x_3)
\psi^{(2)}_{k}(x_2-y_2)\right)\right\|_{\infty}\le
\frac{A}{|R^{\#}||J_R|^k},
\end{gather*}
(where $A$ is a constant only depending on $\psi$)and H\"{o}lder's inequality and similarly for $\widetilde{a_R^{\#}}$. From the definition
of $d_R$ and (\ref{eq3: claim:f}), it is clear that
\begin{align*}
\sum_{R\subset {\R}_i}d_R^2|R||\hat{I_R}|^2 =c\sum_{R\in
{\R}_i}\frac{|R|(\int_R|f_{j,k}(y)|^2dy)}{|R^{\#}|^22^{2i}
|\Omega_i|^{2/p}}\cdot|R|\cdot|\hat{I_R}|^2\le C|\Omega_i|^{1-2/p}.
\end{align*}
Thus we complete the proof of necessity.

Now we turn to the proof of sufficiency part of Theorem \ref{main
flag atom:f}. It suffices to show that for each flag atom $a$,
$$
\|g_F(a)\|_{L^p({\rr}^n\times {\rr}^m)}\le C
$$
where the constant $C$ is independent of $a$. Let $a$ be a flag atom
supported in an open set $\Omega$ with finite measure. Following Chang-Fefferman's idea in [3], set
$$
\overline{\Omega}=\bigcup_{R\subset \Omega \atop R\, {\rm dyadic}}
6R.
$$
We will split the desired estimate into two parts as follows.
\begin{eqnarray}\label{eq1:main part:f}
\int_{\{M(\chi_{\overline{\Omega}})(x)<1/4\}}g_F(a)^p(x)dx\le C,\ \mbox{and}\
\int_{\{M(\chi_{\overline{\Omega}})(x)\ge1/4\}}g_F(a)^p(x)dx\le C.
\end{eqnarray}

We first show the first inequality in (\ref{eq1:main part:f}). Suppose $a=\sum_{R}a_R$ where
$a_R$ are the flag particles. For two different dyadic rectangle
$R=I_R\times J_R$ and $S=I_S\times J_S$, let
$$
m(R,S)=\frac{\min(|I_R|,|I_S|)\,\min(|J_R|,|J_S|)
}{\max(|I_R|,|I_S|)\max(|J_R|,|J_S|)}.
$$

Next we shall give the pointwise estimate for $g_F(a)(x)$. Let $S_x$
 denote the set of all rectangles $S=I_S\times J_S$
containing $x$. Note that for each $x\in {\rr}^n\times {\rr}^m$ and
all $j,k\in {\zz}$, there exists a unique rectangle $S\in S_x$ with
$l(I_S)=2^j$ and $l(J_S)=2^k$. Thus, we may write
$$
g_F^2(a)(x)=\sum_{j,k}|\psi_{j,k}*f(x)|^2=\sum_{S\in
S_x}|\psi_{j,k}*f(x)|^2
$$
where $l(I_S)=2^{j_S}=2^j$ and $l(J_S)=2^{k_S}=2^k$.

Now fix $x$ with $M_s(\chi_{\overline{\Omega}})(x)<1/4$. Then fix
$S=I(S)\times J(S)\in S_x$ and suppose $R=I(R)\times J(R)\subset
\Omega$ is a dyadic rectangle with $\tilde{ R}\cap \tilde{S}\not=
\emptyset$ where $\tilde{R}$ denote the triple of $R$ (note that if
$\tilde{ R}\cap \tilde{S}= \emptyset$, we clearly have
$\psi_{j_S,k_S}*a_R(x)=0$). Then there are four types of such
rectangles $R$.

\medskip

(i) $|I_R|\ge |I_S|$, $|J_R|\ge |J_S|$.

Thus $S\subset 6R$. Hence
$$
1=\frac{|S\cap 6R|}{|S|}\le \frac{|S\cap \overline{\Omega}
|}{|S|}\le M_s(\chi_{\overline{\Omega}})(x)<\frac{1}{4}
$$
which is a contradiction. So this case could not occur.

\medskip

(ii) $|I_R|\le |I_S|$, $|J_R|\ge |J_S|$.

First consider the case $n=m=1$. If $j_R\le k_R$, then for $S\in
S_x$, i.e. $(x_1,x_2)\in I_S\times J_S$, we have
\begin{eqnarray*}
\begin{split}
&|a_R*\psi_{j_S,k_S}(x_1,x_2)|\\
\le& \int |\iiint
a_R^{\#}(y_1,y_3;y_2)\psi_{j_S}^{(1)}(x_1-y_1,x_3-y_3)
\psi_{k_S}^{(2)}(x_2-x_3-y_2)dy_1dy_2dy_3|dx_3.
\end{split}
\end{eqnarray*}
Since
$$
|x_3-c_{I_R}|\le |x_3-y_3|+|y_3-c_{\hat{I_R}}|\le |\hat{I_S}|+
|\hat{I_R}|\le 2|\hat{I_S}|,
$$
 $a_R$ in $x_3$-variable is supported in a translation of  $2\hat{I_S}$. For simplicity, we still use $2\hat{I_S}$ to denote the
 support. Then we have
\begin{eqnarray*}
\begin{split}
&|a_R*\psi_{j_S,k_S}(x_1,x_2)|\\=& \int_{6\hat{I_S}}
|\iiint_{6(I_R\times \hat{I_R}\times J_S)}
(a_R^{\#}(y_1,y_3;y_2)-P(y_1,y_3,c_{J_S}-x_3))\\
&\times
(\psi_{j_S}^{(1)}(x_1-y_1,x_3-y_3)-Q(x_1-c_{I_R},x_3-c_{\hat{I_R}}))
\psi_{k_S}^{(2)}(x_2-x_3-y_2)dy_2dy_3dy_1|dx_3\\
\lesssim & \left(\frac{|J_S|^{k+1}d_R}{|J_R|^{k+1}}\cdot
\frac{|I_R|^{k+1}}{|I_S||\hat{I_S}||I_S|^{k+1}}\cdot
\frac{1}{|J_S|}\right)\,|I_R|\,|\hat{I_R}|\,|J_S|\,|\hat{I_S}|\\
=&\left(\frac{|J_S|}{|J_R|}\cdot
\frac{|I_R|}{|I_S|}\right)^{k+1}\,(d_R|\hat{I_R}|)\,
\left(\frac{|I_R|}{|I_S|}\right),
\end{split}
\end{eqnarray*}
where $P_k(y_1,y_3,c_{J_S}-x_3)$ is the $k^{th}$ order Taylor
polynomial of ${a_R^{\#}}(y_1,y_3,y_2)$ in the third variable at
$(y_1,y_3,c_{J_S}-x_3)$, i.e.
$$
P_k(y_1,y_3,c_{J_S}-x_3)=\sum_{l=0}^k
(y_2+x_2-c_{J_S})^l\frac{\partial^l}{\partial
y_2^l}a_R^{\#}(y_1,y_3,c_{J_S}-x_3)
$$
and $Q_k^{(1)}(x_1-c_{I_R},x_3-c_{\hat{I_R}})$ is $k$-th order
Taylor's polynomial of $\psi_{j_S}^{(1)}(x_1-y_1,x_3-y_3)$ at
$(x_1-c_{I_R},x_3-c_{\hat{I_R}})$.

If $j_R>k_R$. Then we write
\begin{align*}
\begin{split}
&|a_R*\psi_{j_S,k_S}(x_1,x_2)|\\
\le &\int|\iiint
\widetilde{a_R^{\#}}(y_1,y_2;y_3)\psi_{j_S}^{(1)}(x_1-y_1,x_2-x_3-y_2)
\psi_{k_S}^{(2)}(x_3-y_3)dy_1dy_2dy_3|dx_3.\\
\end{split}
\end{align*}
Since $|J_S|\le |J_R|\le |I_R|\le |I_S|$, we have
$$
|x_3-c_{J_R}|\le |x_3-y_3|+|y_3-c_{J_R}|\le 2|\hat{I_S}|.
$$
Hence
\begin{align*}
\begin{split}
|a_R*\psi_{j_S,k_S}(x_1,x_2)| \le
&\int_{6\hat{I_S}}|\iiint_{6(I_R\times \hat{I_R}\times J_S)}
[\widetilde{a_R^{\#}}(y_1,y_2;y_3)-\widetilde{P_k}(y_1,y_2;x_3)]\\
&\times
[\psi_{j_S}^{(1)}(x_1-y_1,x_2-x_3-y_2)-Q_k^{(1)}(x_1-y_1,x_2-x_3-c_{\hat{I_R}})]\\
&\times
\psi_{k_S}^{(2)}(x_3-y_3)dy_1dy_2dy_3|dx_3\\
\le &\left(\frac{|J_S|}{|J_R|}\cdot
\frac{|I_R|}{|I_S|}\right)^{k+1}\,(d_R|\hat{I_R}|)\,\left(\frac{|I_R|}{|I_S|}\right),
\end{split}
\end{align*}
where $\widetilde{P_k}$ denotes the $k$-th order Taylor polynomial
of $\widetilde{a_R^{\#}}$.

Notice that in this special case (ii),
$$
\frac{|J_S|}{|J_R|}\cdot \frac{|I_R|}{|I_S|}=m(S,R)
$$
and
$$
\frac{|I_R|}{|I_S|}=\frac{|cI\cap
cS|}{|cS|}=\frac{1}{|cS|}\int_{cS}\chi_{cR}(t)dt\le M(\chi_{cR})(x)
$$
for some constant $c>0$.

Hence for any $r>0$, we have
\begin{align*}
\begin{split}
|a_R*\psi_{j_S,k_S}(x_1,x_2)|&\lesssim
m(S,R)^{k+1}(d_R|\hat{I_R}|)\left(\frac{|I_R|}{|I_S|}\right)^{1/r}
\left(\frac{|I_R|}{|I_S|}\right)^{1-1/r}\\
&\le m(S,R)^{k+1}(d_R|\hat{I_R}|)M^{1/r}(\chi_{cR})(x)\cdot
\left(\frac{|I_R|}{|I_S|}\right)^{1-1/r}.\\
\end{split}
\end{align*}

In the general case, the moment conditions of $ a_R^{\#}$ allow  us
to use Taylor's formula of higher order to obtain better estimates. We
omit the details.

\medskip

(iii) $|I_R|\ge |I_S|$, $|J_R|\le |J_S|$.

Similarly as in the case (ii), we have
$$
|a_R*\psi_{j_S,k_S}(x_1,x_2)| \le m(S,R)^{k+1}\,(d_R|\hat{I_R}|)\,
M^{1/r}(\chi_{cR})(x)\cdot \left(\frac{|J_R|}{|J_S|}\right)^{1-1/r}.
$$

\medskip
(iv) $|I_R|\le |I_S|$, $|J_R|\le |J_S|$.

We only give the estimate for $n=m=1$. We consider two subcases.

If $j_R<k_R$, then
\begin{align*}
\begin{split}
&|a_R*\psi_{j_S,k_S}(x_1,x_2)|\\
\le& \int |\iiint
a_R^{\#}(y_1,y_3;y_2)\psi_{j_S}^{(1)}(x_1-y_1,x_3-y_3)
\psi_{k_S}^{(2)}(x_2-x_3-y_2)dy_1dy_2dy_3|dx_3.\\
\end{split}
\end{align*}
Since
$
|x_3-c_{\hat{I_R}}|\le |x_3-y_3|+|y_3-c_{\hat{I_R}}|\le
|I_S|+|I_R|\le 2|I_S|,
$ applying Taylor's formula,
we see that this is
\begin{align*}
\begin{split}
\lesssim& \left(\frac{|I_R|^{k+1}}{|I_S||\hat{I_S}||I_S|^{k+1}}\cdot
\frac{|J_R|^{k+1}}{|J_S||J_S|^{k+1}}\cdot
|\hat{I_S}|\right)\left(\iiint|a_R^{\#}(y_1,y_3,y_2)|dy_1dy_2dy_3\right)\\
\le &\left(\frac{1}{|I_S||\hat{I_S}||J_S|}\iiint
|a_R^{\#}(y)|^tdy\right)^{1/t}\left(\frac{1}{|I_S||\hat{I_S}||J_S|}\iiint
\chi_{6(I_R\times \hat{I_R}\times
J_R)}^r(y)dy\right)^{1/r}\\
&\times |\hat{I_S}|\cdot
\left(\frac{|I_R||J_R|}{|I_S||J_S|}\right)^{k+1}\\
\lesssim & d_R\left(\frac{|I_R||\hat{I_R}||J(R)|}{|I_S|
|\hat{I_S}||J_S|}\right)^{1/t}\left(\frac{|\hat{I_R}|}{|\hat{I_S}|}
\right)^{1/r} M^{1/r}(\chi_R)(x)|\hat{I_S}|\cdot
\left(\frac{|I_R||J_R|}{|I_S||J_S|}\right)^{k+1}\\
=&(d_R|\hat{I_R}|)\,M^{1/r}(\chi_{6R})(x)\,m(S,R)^{k+2-1/r}
\end{split}
\end{align*}
for any $r> 1$ and $1/t+1/r=1$.

If $j_R\ge k_R$, then for any $r,t$ with $r> 1$ and $1/t+1/r=1$, we have the similar estimate
\begin{align*}
\begin{split}
&|a_R*\psi_{j_S,k_S}(x_1,x_2)|\\
\lesssim& \left(\frac{|I_R|^{k+1}}{|I_S||\hat{I_S}||I_S|^{k+1}}\cdot
\frac{|J_R|^{k+1}}{|J_S||J_S|^{k+1}}\cdot
|\hat{I_S}|\right)\left(\iiint|\widetilde{a_R^{\#}}(y_1,y_3,y_2)|dy_1dy_2dy_3\right)\\
\le &\left(\frac{1}{|I_S||\hat{I_S}||J_S|}\iiint
|\widetilde{a_R^{\#}}(y)|^tdy\right)^{1/t}\left(\frac{1}{|I_S||\hat{I_S}||J_S|}\iiint
\chi_{6(I_R\times \hat{I_R}\times
J_R)}^r(y)dy\right)^{1/r}\\
&\times |\hat{I_S}|\cdot
\left(\frac{|I_R||J_R|}{|I_S||J_S|}\right)^{k+1}\\
\lesssim &
d_R\left(\frac{|I_R||\hat{I_R}||J(R)|}{|I_S||\hat{I_S}||J_S|}\right)^{1/t}
\left(\frac{|\hat{I_R}|}{|\hat{I_S}|} \right)^{1/r}
M^{1/r}(\chi_R)(x)|\hat{I_S}|\cdot
\left(\frac{|I_R||J_R|}{|I_S||J_S|}\right)^{k+1}\\
=&(d_R|\hat{I_R}|)\,M^{1/r}(\chi_{6R})(x)\,m(S,R)^{k+2-1/r}.
\end{split}
\end{align*}

\medskip

Combing our estimates (i)--(iv), we see that for any $x\in S$,
\begin{align*}
|\psi_{j_S,k_S}*a(x)|^2
\lesssim &\left(\sum_{R\in
\Omega}m^{k+2-1/r}(S,R)(d_R^2|\hat{I_R}|^2)M^{2/r}(\chi_{cR})(x)\right)\left(\sum_{R\subset
\Omega}\Delta_{R,S}\right),
\end{align*}
where
$$
\Delta_{R,S}= \left\{
\begin{aligned}
&\left(\frac{|I_R|}{|I_S|}\right)^{k+2-1/r}\left(\frac{J_S}{J_R}\right)^{k+1/r}\quad
& R\in {\rm type} \
{\rm(ii)}\\
&\left(\frac{|J_R|}{|J_S|}\right)^{k+2-1/r}\left(\frac{I_S}{I_R}\right)^{k+1/r}\quad
&R \in {\rm type}\
{\rm(iii)}\\
&(m(S,R))^{k+2-1/r}\quad  &R \in {\rm type}\
{\rm(iv)}\\
\end{aligned}
\right. .
$$
Using the geometric argument as in \cite{cf}, we get $\sum_{R\in
\Omega}\Delta_{R,S}\le C$. Hence for each $x$ with
$M_s(\chi_{\overline{\Omega}})(x)<1/4$ and for each $r>1$,
\begin{align*}
\begin{split}
g_F^2(a)(x)
\lesssim & \sum_{S\in S_x}\sum_{R\subset
\Omega}m^{k+2-1/r}(S,R)(d_R^2|\hat{I_R}|^2)M^{2/r}(\chi_{cR})(x)\\
= &  \sum_{R\subset \Omega}\left[\sum_{S\in
S_x}m^{k+2-1/r}(S,R)\right](d_R^2|\hat{I_R}|^2)M^{2/r}(\chi_{cR})(x)\\
\le & \sum_{R\subset \Omega}M^{k+2-1/r}(\chi_{\Omega})(x)
(d_R^2|\hat{I_R}|^2)M^{2/r}(\chi_{cR})(x)\\
=& M^{k+2-1/r}(\chi_{\Omega})(x)\sum_{R\subset \Omega}
(d_R^2|\hat{I_R}|^2)M^{2/r}(\chi_{cR})(x).\\
\end{split}
\end{align*}
Then by H\"{o}lder's inequality, for all $s$ and $t$ with
$1/s+1/t=1$
\begin{align*}
\begin{split}
&\int_{\{M(\chi_{\overline{\Omega}})(x)<1/4\}}g_F^p(a)(x)dx\\
\le &\left(\int
M^{p/2(k+2-1/r)s}(\chi_{\Omega})(x)dx\right)^{1/s}\left(\int(\sum_{R\subset
\Omega}(d_R^2|\hat{I_R}|^2) M^{2/r}(\chi_{R})(x))^{(p/2)\cdot
t}dx\right)^{1/t}.
\end{split}
\end{align*}
 If we choose $t=2/p$, choose $r$ so that $1<r<2$
and finally choose $k$ such that $p/2(k+2-1/r)s>1$ (indeed,
$k=[2/p-3/2]$ would be sufficient). Then
\begin{align*}
\begin{split}
\int_{\{M(\chi_{\overline{\Omega}})(x)<1/4\}}g_F^p(a)(x)dx &\le
\left(\int
M^{p/2(k+2-1/r)s}(\chi_{\Omega})(x)dx\right)^{1/s}(\sum_{R\subset
\Omega}d_R^2|\hat{I_R}|^2|R|)^{1/t}\\
&\le C|\Omega|^{1/s}|\Omega|^{(1-2/p)1/t} \le C
\end{split}
\end{align*}
which gives the first inequality in (\ref{eq1:main part:f}).

As for the second inequality in (\ref{eq1:main part:f}), note that $pt=2$, by H\"{o}lder's
inequality,
\begin{align*}
\begin{split}
\int_{\{M(\chi_{\overline{\Omega}})(x)> 1/4\}}g_F^p(a)(x)dx \le
&\,\left(\int
g_F^{pt}(a)(x)dx\right)^{1/t}|\{M(\chi_{\overline{\Omega}})(x)>1/4\}|^{1/s}\\
\le & \,C(\|a\|_2^{2/t})\,|\Omega|^{1-1/t}\\
\le &\,C\|a\|_2^{p/2}\,|\Omega|^{1-p/2}\\
\le &\,C|\Omega|^{(1-2/p)p/2}|\Omega|^{1-p/2}\le C.
\end{split}
\end{align*}
To finish the proof of Theorem 1.1, let $f$ be in both $L^2$ and $H^p_F(\rrnm).$ Then, as before, $f=\sum_Rf_R=\sum_i \sum_{R\subset {\R}_i}f_R=\sum_i \lambda_i a_i,$ where $a_i$ are $H^p_F(\rrnm)$ atoms. We only need to show that the series of such an atomic decomposition of $f$ converges also in $L^2$ norm. To this end, as above, by the duality argument as in \eqref{dual},
\begin{eqnarray*}
\begin{split}
\|\sum_{|i|>N} \lambda_i a_i\|_{2}=&\sup_{\|g\|_2\le 1}\left|\int_{{\rr}^m\times
{\rr}^n}\sum_{|i|>N} \lambda_i a_i(x,y)g(x,y)dx dy\right|\\
=&\sup_{\|g\|_2\le 1}\left|\int_{{\rr}^m\times
{\rr}^n}\sum_{|i|>N} \sum_{R\in\R_i}f_R(x,y) g(x,y)dx dy\right|\\
\le& C\left\{\sum_{|i|>N}[\sum_{R=I\times J\in \R_i}|R|\cdot\sup_{(u,v)\in
R}|\psi_{j,k}\ast f(u,v)|]^2 \right\}^{1/2}.
\end{split}
\end{eqnarray*}
Applying the estimate in (2.3) yields
$$
\|\sum_{|i|>N} \lambda_i a_i\|_{2}\le C
\sum_{|i|>N} 2^{2i}|\Omega_i|,$$
where the last term above tends to zero as $N$ goes to infinity since $f\in L^2.$ $\qe$

The proof of Theorem \ref{main flag atom:f} is concluded.

\section{Proof of Theorem \ref{main:application}}
 We first assume that $f\in L^2\cap H^p_F(\rrnm)$. Let  $f=\sum_{i=1}^\infty\lambda_ia_i$ be an atomic decomposition of $f,$ where the series converges in both $L^2$ and $H^p_F.$
Since $T$ is bounded on $L^2,$ thus,
$$\left\|T[\sum_{i=1}^\infty\lambda_ia_i]\right\|^p_{L^p}=\left\|\sum_{i=1}^\infty\lambda_i (Ta_i)\right\|^p_{L^p}\lesssim \sum_{i=1}^\infty\lambda_i^p\|Ta_i\|^p_p.$$
This implies that if $\|Ta\|_p\le C$ for all $H^p_F$ atoms $a$ then $\|Tf\|_p\le C\|f\|_{H^p_F}$ for $f\in L^2\cap H^p_F.$ Similarly, if $\|Ta\|_{H^p_F}\le C$ for all $H^p_F$ atoms $a$ then  $\|Tf\|_{H^p_F}\le C\|f\|_{H^p_F}$ for $f\in L^2\cap H^p_F.$ Since $L^2\cap H^p_F$ is dense in $H^p_F$, a density argument finishes the proof of sufficiency.
The necessity is obvious because $\|a\|_{H^p_F}\lesssim 1.$
$\qe$

\end{document}